\documentclass[11pt,reqno]{amsart}
\usepackage{enumerate}
\usepackage{amsfonts, amsmath, amssymb, amscd, amsthm, bm}
\usepackage{enumerate}
\usepackage{url}
\usepackage[linktocpage=true,colorlinks,citecolor=magenta,linkcolor=blue,urlcolor=magenta]{hyperref}
\usepackage{multicol}
\usepackage{comment}

 \newtheorem{thm}{Theorem}[section]
 \newtheorem{cor}[thm]{Corollary}
  
 \newtheorem{lem}[thm]{Lemma}
 \newtheorem{assump}[thm]{Assumption}
 
 \theoremstyle{definition}

 \newtheorem{rem}[thm]{Remark}
 \newtheorem*{ack}{Acknowledgments}

 \theoremstyle{claim}
 \numberwithin{equation}{section}
\begin{document}
\title[New pinching estimates for Inverse curvature flows]{New pinching estimates for Inverse curvature flows in space forms}
\author[Y. Wei]{Yong Wei}
\address{Mathematical Sciences Institute,
Australian National University, Canberra,
ACT 2601 Australia}
\email{\href{mailto:yong.wei@anu.edu.au}{yong.wei@anu.edu.au}}

%


\begin{abstract}
We consider the inverse curvature flow of strictly convex hypersurfaces in the space form $N$ of constant sectional curvature $K_N$ with speed given by $F^{-\alpha}$, where $\alpha\in (0,1]$ for $K_N=0,-1$ and $\alpha=1$ for $K_N=1$,  $F$ is a smooth, symmetric homogeneous of degree one function which is inverse concave and has dual $F_*$ approaching zero on the boundary of the positive cone $\Gamma_+$. We show that the ratio of the largest principal curvature to the smallest principal curvature of the flow hypersurface is controlled by its initial value. This can be used to prove the smooth convergence of the flows.
\end{abstract}

\maketitle

\section{Introduction}

Let $N$ be a real simply connected space form of dimension $n+1$ with constant sectional curvature $K_N$, i.e., $N$ is the Euclidean space $\mathbb{R}^{n+1}$ if $K_N=0$, $N$ is the hyperbolic space $\mathbb{H}^{n+1}$ if $K_N=-1$ and $N$ is the sphere $\mathbb{S}^{n+1}$ if $K_N=1$. Let $M_0$ be a smooth closed and strictly convex hypersurface of dimension $n\geq 2$ in $N$ parameterized by the embedding $X_0:M^n\to X_0(M^n)=M_0\subset N$. We consider the inverse curvature flow starting from $M_0$, which is a smooth family of embeddings $X:M^n\times[0,T)\to N$ satisfying
\begin{equation}\label{s1:flow1}
 \left\{\begin{aligned}
 \frac{\partial}{\partial t}X(x,t)=&~F(\mathcal{W}(x,t))^{-\alpha}\nu(x,t),\\
 X(\cdot,0)=&~X_0(\cdot),
  \end{aligned}\right.
 \end{equation}
 where $\alpha>0$, $\mathcal{W}(x,t)$ is the matrix of the Weingarten map of $M_t=X(M^n,t)$, $F(\mathcal{W})=f(\kappa(\mathcal{W}))$ is a smooth symmetric function of the principal curvatures $\kappa=(\kappa_1,\cdots,\kappa_n)$ of $M_t$, and $\nu$ is the outer unit normal of $M_t$.

 We will use the following assumptions on the curvature function $F$.

 \begin{assump}\label{assum-1}
Let $\Gamma_+=\{(\kappa_1,\cdots,\kappa_n)\in \mathbb{R}^n: \kappa_i>0, i=1,\cdots,n\}$ be the
positive quadrant in $\mathbb{R}^n$. Assume that
\begin{itemize}
  \item[(i)] $F(\mathcal{W})=f(\kappa(\mathcal{W}))$, where $\kappa(\mathcal{W})$ gives the eigenvalues of $\mathcal{W}$ and $f$ is a smooth symmetric function on $\Gamma_+$.
  \item[(ii)] $f$ is strictly increasing in each argument, i.e., ${\partial f}/{\partial \kappa_i}>0$ on $\Gamma_+$, $\forall~i=1,\cdots,n$.
  \item[(iii)] $f$ is homogeneous of degree $1$, i.e., $f(k\kappa)=kf(\kappa)$ for any $k>0$.
  \item[(iv)] $f$ is strictly positive on $\Gamma_+$ and is normalized such that $f(1,\cdots,1)=n$.
  \item[(v)] Either:
  \begin{itemize}
    \item[(a)] $f$ is concave and $f$ approaches zero on the boundary of $\Gamma_+$.
    \item[(b)]  $f$ is concave and inverse concave, i.e., the dual function
    \begin{equation}\label{s1:f-dual}
      f_*(x_1,\cdots,x_n)=f(x_1^{-1},\cdots,x_n^{-1})^{-1}
    \end{equation}
    is also concave.
    \item[(c)] $f$ is inverse concave and $f_*$ approaches zero on the boundary of $\Gamma_+$.
 \end{itemize}
\end{itemize}
\end{assump}

Before we state our main result, let's review the known results on the behavior of the inverse curvature flow \eqref{s1:flow1} in Euclidean space with curvature function $f$ satisfying (i)-(iv) and different items in (v) of Assumption \ref{assum-1}. If the initial hypersurface $M_0$ is strictly convex, then for $\alpha\in (0,1]$ and $f$ satisfies either (b) or (c), Urbas \cite{Urbas91JDG} proved that the flow \eqref{s1:flow1} has a smooth strictly convex solution $M_t$ for all $t\in [0,\infty)$, $M_t$ expands to the infinity and properly rescaled flow hypersurfaces converge smoothly to a round sphere. For $\alpha\in (0,1]$ and the dimension $n=2$, the same result was proved by the author with Li and Wang in \cite{LiWangWei2016} without any second derivative assumption on $f$. Note that the case $n=2$ was also considered by Urbas in \cite{Urbas91JDG} provided that $f$ is inverse concave.  The key ingredients in \cite{LiWangWei2016} are the pinching estimate on the principal curvatures of evolving hypersurfaces $M_t$ and Andrews' \cite{And04} regularity estimate for parabolic equations in two variables without any concavity assumption.  For $\alpha>1$, the case (a) was considered by Gerhardt \cite{Gerhardt2014} and the case that $f$ is concave and the initial hypersurface $M_0$ is sufficiently pinched was considered by Kr\"{o}ner-Scheuer \cite{Kron-Sche2017}.

For star-shaped (not necessarily convex) hypersurface $M_0$, Gerhardt \cite{Gerhardt1990,Gerhardt2014} ($\alpha\in (0,1]$) and Urbas \cite{urbas1990expansion} ($\alpha=1$) studied the flow \eqref{s1:flow1} with concave curvature function $f$ which satisfies $f|_{\Gamma}>0$ and $f|_{\partial\Gamma}=0$ for an open convex symmetric cone $\Gamma$ containing the positive cone $\Gamma_+$, and proved a similar convergence result. Scheuer \cite{Sch2016} improved the asymptotical behavior of the flow \eqref{s1:flow1} considered in \cite{Gerhardt2014} by showing that the flow becomes close to a flow of a sphere.

The inverse curvature flow has also been studied in other ambient spaces, in particular in the hyperbolic space and in sphere, with curvature function $f$ satisfying (i)-(iv) and item (a) or (b) in (v) of Assumption \ref{assum-1}.  See \cite{Ding-11,gerhardt2011-H^n,Gerhardt2015,LiWangWei2016,makowski2013rigidity,scheuer2015-gradient,Scheuer2015-1}.

The main purpose of this paper is to show the following pinching estimate along the flow \eqref{s1:flow1} with $F$ satisfying (i)--(iv) and (v)(c) in Assumption \ref{assum-1}.
 \begin{thm}\label{s1:thm-pinc}
Let $M_0$ be a smooth, closed and strictly convex hypersurface of dimension $n\geq 2$ in the simply connected space form $N$ of constant sectional curvature $K_N$.  Assume that $F$ satisfies (i)--(iv) and (v)(c) in Assumption \ref{assum-1}. Then along the flow \eqref{s1:flow1} with $\alpha\in(0,1]$ for $K_N=0,-1$ and $\alpha=1$ for $K_N=1$, the principal curvatures $\kappa=(\kappa_1,\cdots,\kappa_n)$ of $M_t$ satisfies
 \begin{equation}\label{s1:pinc2}
   \kappa_{n}\leq C\kappa_{1}
 \end{equation}
for all $t\in[0,T)$, where $\kappa_i$ are labelled such that $\kappa_1\leq \kappa_2\leq\cdots\leq\kappa_n$ and $C>0$ is a constant depending only on the initial data $M_0$.
\end{thm}

The proof of Theorem \ref{s1:thm-pinc} will be given in section \ref{sec:pinc}. By applying Andrews' tensor maximum principle \cite[\S 3]{Andrews2007},  we firstly prove that if $F$ is inverse concave, then the minimum of the smallest eigenvalue of $F^{-1}h_i^j$ over $M_t$ is non-decreasing in time $t$ along the flow \eqref{s1:flow1} with $\alpha\in(0,1]$ for $K_N=0,-1$ and $\alpha=1$ for $K_N=1$. Combining this with the assumption that $f_*$ approaches zero on the boundary of the positive cone $\Gamma_+$ yields the pinching estimate \eqref{s1:pinc2}. We remark that Kr\"{o}ner-Scheuer \cite{Kron-Sche2017} recently constructed an example to show that the strictly convexity (in particular the pinching estimate) may be lost along the flow \eqref{s1:flow1} in Euclidean space with $\alpha>1$ and $F=H$. Thus the range $\alpha\in (0,1]$  is the best to hope for when it comes to the pinching
estimates.

The pinching estimate \eqref{s1:pinc2} can be used to prove the convergence result of the flow \eqref{s1:flow1}. As we mentioned above, Urbas \cite{Urbas91JDG} has proved the smooth convergence of the flow \eqref{s1:flow1} for strictly convex initial hypersurface with curvature function $F$ satisfying (i)--(iv) and (v)(c) in Assumption \ref{assum-1}. In the following, we state the convergence results for the flow in the hyperbolic space and in sphere.
 \begin{cor}\label{s1:thm-conv}
 Let $M_0$ be a smooth closed and strictly convex hypersurface in $N^{n+1} (n\geq 2)$, where $N$ is either the hyperbolic space $\mathbb{H}^{n+1}$ or the sphere $\mathbb{S}^{n+1}$. Assume that $F$ satisfies (i)--(iv) and (v)(c) in Assumption \ref{assum-1}. Then,
\begin{itemize}
  \item[(i)]  along the flow \eqref{s1:flow1} in $\mathbb{H}^{n+1}$ with $\alpha\in(0,1]$, the solution $M_t$ exists for all time $t\in[0,\infty)$, $M_t$ expands to infinity as $t\to \infty$ and properly rescaled hypersurfaces converge smoothly to a geodesic sphere.
  \item[(ii)] along the flow \eqref{s1:flow1} in $\mathbb{S}^{n+1}$ with $\alpha=1$, the solution $M_t$ exists for $t\in[0,T)$ with $T<\infty$, $M_t$ expands to the equator of the sphere $\mathbb{S}^{n+1}$ smoothly as $t\to T$. Moreover, if $F_*$ is strictly concave or $F_*$ equals to the mean curvature $H$, a properly rescaled hypersurfaces converge smoothly to a geodesic sphere.
\end{itemize}
\end{cor}
The key step to prove Corollary \ref{s1:thm-conv} is to derive the $C^{2,\alpha}$ estimate of the evolving hypersurfaces. Most of the previous papers (see \cite{Gerhardt1990,gerhardt2011-H^n,Gerhardt2014,Gerhardt2015,Kron-Sche2017,Scheuer2015-1}) required the concavity of the speed function so that the $C^{2,\alpha}$ estimate by Krylov and Evans can be applied. In our case, we can replace the concavity $F$ by the inverse concavity of $F$. The idea can be described easily in the Euclidean case (as was used by Urbas \cite{Urbas91JDG}): Since the evolving hypersurfaces $M_t$ are strictly convex, we can reparametrise them using the Gauss map and write the flow \eqref{s1:flow1} as a scalar parabolic PDE for the support function $s(z,t)$ of $M_t$
\begin{equation}\label{s1:flow-s}
  \frac{\partial }{\partial t}s(z,t)~=~F_*^{\alpha}(\bar{\nabla}_i\bar{\nabla}_js(z,t)+s(z,t)\bar{g}_{ij}), \quad (z,t)\in \mathbb{S}^n\times [0,T),
\end{equation}
where $\bar{g}_{ij}$ and $\bar{\nabla}$ are round metric and its Levi-Civita connection on $\mathbb{S}^n$. Using the inverse concavity of $F$, the right hand side of \eqref{s1:flow-s} is now concave with respect to the second spatial derivatives of $s$ provided that $\alpha\leq 1$. The theorem of Krylov and Evans can now be applied to derive the $C^{2,\alpha}$ regularity. This idea was also used in \cite{And2000,Andrews-McCoy-Zheng} to study the contracting curvature flow of convex hypersurfaces in the Euclidean space.  In the hyperbolic space and in sphere, we can also develop the similar idea. In \cite{And-Wei2017}, the author with Andrews derived the $C^{2,\alpha}$ regularity of the volume preserving flow in the hyperbolic space by writing the flow as a scalar parabolic PDE of the support function of the corresponding evolution in the unit Euclidean ball, which is concave with respect to the second spatial derivatives.  To prove Corollary \ref{s1:thm-conv}, we will employ the same formulation to derive the $C^{2,\alpha}$ regularity estimate. The flow in the sphere can also be treated in a similar way.

\begin{rem}
Though the convergence result of the flow \eqref{s1:flow1} in the Euclidean space for curvature functions satisfying (i)--(iv) and (v)(c) in Assumption \ref{assum-1}  has been proved by Urbas in \cite{Urbas91JDG}, the pinching estimate \eqref{s1:pinc2} given in Theorem \ref{s1:thm-pinc} is still interesting. The  convergence result in \cite{Urbas91JDG} was proved by considering the evolution of the normalized support function $\tilde{s}=e^{-t}s$ of the evolving hypersurfaces
\begin{equation}\label{s1:flow-s}
  \frac{\partial \tilde{s}}{\partial t}~=~F_*^{\alpha}(\bar{\nabla}_i\bar{\nabla}_j\tilde{s}+\tilde{s}\bar{g}_{ij})-\tilde{s}
\end{equation}
on $\mathbb{S}^n\times [0,\infty)$ and derived the $C^2$ estimate of $\tilde{s}$ using the maximum principle directly. The pinching estimate \eqref{s1:pinc2} together with the bounds on the speed function $F$ (which are immediate consequences of the maximum principle) provides an alternative way to derive the $C^2$ estimate of the flow.
\end{rem}

\begin{ack}
The author would like to thank the referees for carefully reading of this manuscript and providing many helpful suggestions.  The author was supported by Ben Andrews throughout his Australian Laureate Fellowship FL150100126 of the Australian Research Council.
\end{ack}

\section{Preliminaries}\label{sec:pre}

For a smooth symmetric function $F(A)=f(\kappa(A))$, where $A=(A_{ij})\in \mathrm{Sym}(n)$ is a symmetric matrix and $\kappa(A)=(\kappa_1,\cdots,\kappa_n)$ give the eigenvalues of $A$, we denote by $\dot{F}^{ij}$ and $\ddot{F}^{ij,kl}$ the first and second derivatives of $F$ with respect to the components of its argument, so that
\begin{equation*}
  \frac{\partial}{\partial s}F(A+sB)\bigg|_{s=0}=\dot{F}^{ij}(A)B_{ij}
\end{equation*}
and
\begin{equation*}
  \frac{\partial^2}{\partial s^2}F(A+sB)\bigg|_{s=0}=\ddot{F}^{ij,kl}(A)B_{ij}B_{kl}
\end{equation*}
for any two symmetric matrixs $A,B$.  We also use the notation
\begin{equation*}
  \dot{f}^i(\kappa)=\frac{\partial f}{\partial \kappa_i}(\kappa),\quad  \ddot{f}^{ij}(\kappa)=\frac{\partial^2 f}{\partial \kappa_i\partial\kappa_j}(\kappa)
\end{equation*}
for the derivatives of $f$ with respect to $\kappa$. At any diagonal $A$ with distinct eigenvalues, the second derivative $\ddot{F}$ of $F$ in direction $B\in \mathrm{Sym}(n)$ is given in terms of $\dot{f}$ and $\ddot{f}$ by  (see \cite{And1994-2,Andrews2007}):
\begin{equation}\label{s2:F-ddt}
  \ddot{F}^{ij,kl}B_{ij}B_{kl}=\sum_{i,k}\ddot{f}^{ik}B_{ii}B_{kk}+2\sum_{i>k}\frac{\dot{f}^i-\dot{f}^k}{\kappa_i-\kappa_k}B_{ik}^2.
\end{equation}
This formula makes sense as a limit in the case of any repeated values of $\kappa_i$.

\subsection{Inverse concave function}
For any positive definite symmetric matrix $A\in \mathrm{Sym}(n)$, define $F_*(A)=F(A^{-1})^{-1}$. Then $F_*(A)=f_*(\kappa(A))$, where $f_*$ is defined in \eqref{s1:f-dual}. Since $f$ is defined on the positive definite cone $\Gamma_+$, the following lemma characterizes the inverse concavity of $f$ and $F$.
\begin{lem}[\cite{Andrews2007,And-Wei2017}]
\begin{itemize}
\item[(i)] $f$ is inverse concave if and only if the following matrix
 \begin{equation}\label{s1:matrx}
   \left(\ddot{f}^{ij}+\frac 2{x_i}\dot{f}^i\delta_{ij}\right)\geq~0.
 \end{equation}
  \item[(ii)] $F_*$ is concave if and only if $f_*$ is concave.
   \item[(iii)]If $f$ is inverse concave, then
 \begin{equation}\label{s2:inv-conc-1}
   \sum_{k,l=1}^n\ddot{f}^{kl}y_ky_l+2\sum_{k=1}^n\frac {\dot{f}^k}{\kappa_k}y_k^2~\geq ~2f^{-1}(\sum_{k=1}^n\dot{f}^ky_k)^2
 \end{equation}
 for any $y=(y_1,\cdots,y_n)\in \mathbb{R}^n$, and
\begin{equation}\label{s2:inv-conc}
\frac{\dot{f}^k-\dot{f}^l}{\kappa_k-\kappa_l}+\frac{\dot{f}^k}{\kappa_l}+\frac{\dot{f}^l}{\kappa_k}\geq~0,\quad \forall~k\neq l.
\end{equation}
%
  \item[(iv)] If $f$ is inverse concave, then
  \begin{equation}\label{s2:ic-2}
  \sum_{k=1}^n \dot{f}^k\kappa_k^2\geq~ f^2/n.
  \end{equation}
\end{itemize}
\end{lem}

We can easily see that a convex function $f:\Gamma_+\to\mathbb{R}$ satisfies the case (v)(c) in Assumption \ref{assum-1}. Firstly, the inequality \eqref{s1:matrx} is obviously true since $f$ is convex and strictly increasing.  Secondly, the convexity of $f$ also implies that
\begin{equation*}
  f(x_1,\cdots,x_n)~\geq~\sum_{i=1}^nx_i.
\end{equation*}
This implies  that
\begin{equation*}
  f_*(z_1,\cdots,z_n)=f(\frac 1{z_1},\cdots,\frac 1{z_n})^{-1}\leq \left(\frac 1{z_1}+\cdots+\frac 1{z_n}\right)^{-1}=\frac{E_n(z)}{nE_{n-1}(z)},
\end{equation*}
where $E_n(z)$ and $E_{n-1}(z)$ are the normalized elementary symmetric polynomial of $z=(z_1,\cdots,z_n)$. Thus $f_*$ approaches zero on the boundary of $\Gamma_+$. Some important examples of function satisfying (i)-(iv) and (v)(c) of Assumption \ref{assum-1} are $f=nE_k^{1/k}, k=1,\cdots,n$, and the the power means $H_r=n^{1-\frac 1r}\left(\sum_{i}\kappa_i^{r}\right)^{1/r}, r>0$. More examples can be constructed as follows:  If $G_1$ is homogeneous of degree one, increasing in each argument, and inverse-concave, and $G_2$ satisfies (i)-(iv) and (v)(c) of Assumption \ref{assum-1}, then $F=G_1^\sigma G_2^{1-\sigma}$ satisfies (i)-(iv) and (v)(c) of Assumption \ref{assum-1} for any $0<\sigma<1$  (see \cite{Andrews2007,ALM14} for more examples on inverse concave, or convex functions).

 \subsection{Evolution equations on the hypersurface}\label{sec:2-1}
Set $\Phi(r)=-r^{-\alpha}$ for $r>0$ and denote $\Phi'=\frac d{dr}\Phi(r)$. Then the flow \eqref{s1:flow1} is equivalent to \begin{equation}\label{s2:flow1}
\frac{\partial}{\partial t}X(x,t)=~-\Phi(x,t)\nu(x,t),
\end{equation}
where $\Phi=-F^{-\alpha}$. Let $\nu, g_{ij}$ and $h_{ij}$ be the unit normal vector field, the induced metric and the second fundamental form of the evolving hypersurface $M_t$.  Throughout this paper we will always evaluate the derivatives of $F$ at the Weingarten matrix $\mathcal{W}=(h_i^j)$ and the derivatives of $f$ at the principal curvatures $\kappa$. We have the following evolution equations for $g_{ij}$,  $h_{ij}$ and the speed function $\Phi$ along the flow \eqref{s2:flow1} (see \cite{And1994-3,Gerhardt2006}):
\begin{equation}\label{s2:g-evl}
  \frac{\partial}{\partial t}g_{ij}=~-2\Phi h_{ij},
\end{equation}
\begin{equation}\label{s2:speed}
  \frac{\partial}{\partial t}\Phi=~\Phi'\dot{F}^{ij}\nabla_i\nabla_j\Phi+\Phi\Phi'\left(\dot{F}^{ij}h_i^kh_{kj}+K_N\dot{F}^{ij}g_{ij}\right),
\end{equation}
and
\begin{equation}\label{s2:h}
\begin{aligned}
\frac{\partial}{\partial t}h_{ij}=&~{\Phi}'\dot{F}^{kl}\nabla_k\nabla_lh_{ij}+{\Phi}'\ddot{F}^{kl,mn}\nabla_ih_{kl}\nabla_jh_{mn}+{\Phi}''\nabla_iF\nabla_jF\\
&\quad +{\Phi}'\left(\dot{F}^{kl}h^p_kh_{pl}-K_N\dot{F}^{kl}g_{kl}\right)h_{ij}\\
&\quad-(\Phi+\Phi'F)h^k_ih_{kj}+K_N(\Phi+\Phi'F)g_{ij},
\end{aligned}
\end{equation}
where $\nabla$ denotes the Levi-Civita connection  on $M_t$ with respect to the induced metric, $\Phi'=\alpha F^{-\alpha-1}$ and $\Phi''=-\alpha(\alpha+1)F^{-\alpha-2}$.


\section{Pinching estimate}\label{sec:pinc}
In this section, we prove that the pinching ratio of the largest principal curvature and the smallest principal curvature is controlled by its initial value along the flow \eqref{s1:flow1} with power $\alpha\in (0,1]$ for $K_N=0,-1$ and $\alpha=1$ for $K_N=1$ if the function $f$ is inverse concave and has dual $f_*$ approaching zero on the boundary of $\Gamma_+$.  A similar pinching estimate was proved by Andrews-McCoy-Zheng in \cite[Lemma 11]{Andrews-McCoy-Zheng} for contracting curvature flow
\begin{equation}\label{s3:flow-contr}
\frac{\partial}{\partial t}X(x,t)=~-F(\mathcal{W}(x,t))^{\alpha}\nu(x,t)
\end{equation}
in the Euclidean space $\mathbb{R}^{n+1}$ with power $\alpha=1$ and inverse concave curvature function $f$ with dual $f_*$ approaching zero on the boundary of $\Gamma_+$. Their proof relies on the Gauss map parametrization of the flow \eqref{s3:flow-contr}. A similar computation also appeared in \cite{Mcc2017} for the mixed volume preserving flow in Euclidean space. Our proof of the pinching estimate along the inverse curvature flow is inspired by their argument but instead of using the Gauss map parametrization of the flow we will prove this estimate directly using the evolution equations \eqref{s2:speed} and \eqref{s2:h}. This makes it possible to deal with the flow \eqref{s1:flow1} in hyperbolic space and in sphere.

The main tool to prove the pinching estimate is the tensor maximum principle. For the convenience of readers, we include here the statement of the tensor maximum principle, which was first proved by Hamilton \cite{Ha1982} and was generalized by Andrews \cite{Andrews2007}.
\begin{thm}[\cite{Andrews2007}]\label{s2:tensor-mp}
Let $S_{ij}$ be a smooth time-varying symmetric tensor field on a compact manifold $M$, satisfying
\begin{equation}
\frac{\partial}{\partial t}S_{ij}=a^{kl}\nabla_k\nabla_lS_{ij}+u^k\nabla_kS_{ij}+N_{ij},
\end{equation}
where $a^{kl}$ and $u$ are smooth, $\nabla$ is a (possibly time-dependent) smooth symmetric connection, and $a^{kl}$ is positive definite everywhere. Suppose that
\begin{equation}\label{s2:TM2}
  N_{ij}v^iv^j+\sup_{\Lambda}2a^{kl}\left(2\Lambda_k^p\nabla_lS_{ip}v^i-\Lambda_k^p\Lambda_l^qS_{pq}\right)\geq 0
\end{equation}
whenever $S_{ij}\geq 0$ and $S_{ij}v^j=0$. If $S_{ij}$ is positive definite everywhere on $M$ at $t=0$ and on $\partial M$ for $0\leq t\leq T$, then it is positive on $M\times[0,T]$.
\end{thm}

\begin{lem}\label{s3:lem-pinc}
Assume that $f$ is inverse concave, $\alpha\in(0,1]$ for $K_N=0,-1$ and $\alpha=1$ and $K_N=1$. Then the minimum of the smallest eigenvalue of $F^{-1}h_i^j$ over $M_t$ is strictly increasing in time $t$ along the flow \eqref{s1:flow1} unless $M_t$ is a totally geodesic sphere.
\end{lem}
\proof
Define the tensor $G_{ij}=F^{-1}h_{ij}-Cg_{ij}$, where $C$ is chosen to make $G_{ij}$ is positive definite initially. Since $\Phi(F)=-F^{-\alpha}$,  equation \eqref{s2:speed} implies that
\begin{align}\label{s3:F-evl}
  \frac{\partial}{\partial t}F=&~\Phi'\dot{F}^{kl}\nabla_k\nabla_lF+\Phi''\dot{F}^{kl}\nabla_kF\nabla_lF\nonumber\\
  &\quad+\Phi\left(\dot{F}^{ij}h_i^kh_{kj}+K_N\dot{F}^{ij}g_{ij}\right).
\end{align}
Combining \eqref{s3:F-evl} with \eqref{s2:g-evl} and \eqref{s2:h},  we have the evolution equation for $G_{ij}$ as follows:
\begin{align}\label{s6:G1}
\frac{\partial}{\partial t}G_{ij}=&\alpha F^{-\alpha-1}\dot{F}^{kl}\nabla_k\nabla_lG_{ij}+2\alpha F^{-\alpha-2}\dot{F}^{kl}\nabla_kF\nabla_lG_{ij}\nonumber\\
&\quad +\alpha F^{-\alpha-2}\ddot{F}^{kl,mn}\nabla_ih_{kl}\nabla_jh_{mn}-\alpha(\alpha+1)F^{-\alpha-3}\nabla_iF\nabla_jF\nonumber\\
&\quad +\alpha(\alpha+1)F^{-\alpha-4}\dot{F}^{kl}\nabla_kF\nabla_lFh_{ij}\nonumber\\
&\quad +(\alpha+1)F^{-\alpha-2}\dot{F}^{kl}h^p_kh_{pl}h_{ij}-(\alpha-1)F^{-\alpha-1} h^k_ih_{kj}\nonumber\\
&\quad -2CF^{-\alpha}h_{ij}+K_N(1-\alpha)F^{-\alpha-1}\left(\dot{F}^{kl}g_{kl}F^{-1}h_{ij}-g_{ij}\right).
\end{align}
We will apply the tensor maximum principle to show that if $G_{ij}\geq 0$ initially then it remains true for later time. Let $(x_0,t_0)$ be the point where $G_{ij}$ has a null vector field $v$, i.e., $G_{ij}v^j=0$ at $(x_0,t_0)$. If we choose normal coordinates at $(x_0,t_0)$ such that the Weingarten matrix is diagonalised with eigenvalues $\kappa=(\kappa_1,\dots,\kappa_n)$ in increasing order, then the null vector $v$ is the eigenvector $e_1$ corresponding to the eigenvalue $\kappa_1$. Let's first look at the zero order terms of \eqref{s6:G1}, i.e., the terms in the last two lines of \eqref{s6:G1} which we denote by $Q_0$.
\begin{align}\label{s3:Q0-1}
Q_0v^iv^j=&~ (\alpha+1)F^{-\alpha-2}\dot{F}^{kl}h^p_kh_{pl}h_{ij}v^iv^j-(\alpha-1)F^{-\alpha-1} h^k_ih_{kj}v^iv^j\nonumber\\
&\quad -2CF^{-\alpha}h_{ij}v^iv^j+K_N(1-\alpha)F^{-\alpha-1}\left(\dot{F}^{kl}g_{kl}F^{-1}h_{ij}-g_{ij}\right)v^iv^j\nonumber\\
=&~ (\alpha+1)F^{-\alpha-2}\dot{F}^{kl}h^p_kh_{pl}\kappa_1-(\alpha-1)F^{-\alpha-1} \kappa_1^2\nonumber\\
&\quad -2CF^{-\alpha}\kappa_1+K_N(1-\alpha)F^{-\alpha-1}\left(\dot{F}^{kl}g_{kl}F^{-1}\kappa_1-1\right).
\end{align}
Since
\begin{equation*}
\dot{F}^{kl}g_{kl}F^{-1}\kappa_1-1=\frac{\sum_k\dot{f}^k(\kappa_1-\kappa_k)}{\sum_k\dot{f}^k\kappa_k}\leq 0,
\end{equation*}
the last term of \eqref{s3:Q0-1} is nonnegative provided that $K_N(1-\alpha)\leq 0$. Note that $G_{ij}v^iv^j=0$ implies that $F^{-1}\kappa_1=C$ at $(x_0,t_0)$. Then if $K_N(1-\alpha)\leq 0$, we have
\begin{align*}
Q_0v^iv^j\geq &~ (\alpha+1)CF^{-\alpha-1}\left(\sum_k\dot{f}^k\kappa_k^2-F\kappa_1\right)\\
=&~ (\alpha+1)CF^{-\alpha-1}\sum_k\dot{f}^k\kappa_k\left(\kappa_k-\kappa_1\right)~\geq~0
\end{align*}
at $(x_0,t_0)$ . Therefore $Q_0v^iv^j\geq 0$ at $(x_0,t_0)$ if $K_N(1-\alpha)\leq 0$.

By continuity we can assume that $h_i^j$ has all eigenvalues distinct at $(x_0,t_0)$ and satisfies $\kappa_1<\kappa_2<\cdots<\kappa_n$. This is possible since for any positive definite symmetric matrix $A$ with $A_{ij}\geq CF(A)\delta_{ij}$ and $A_{ij}v^iv^j=CF(A)|v|^2$ for some $v\neq 0$, there is a sequence of symmetric matrixes $\{A^{(k)}\}$ approaching $A$, satisfying $A^{(k)}_{ij}\geq CF(A^{(k)})\delta_{ij}$ and $A^{(k)}_{ij}v^iv^j=CF(A^{(k)})|v|^2$ and with each $A^{(k)}$ having distinct eigenvalues. Hence it suffice to prove the result in the case where all of $\kappa_i$ are distinct. Since the null vector $v=e_1$ and $C=F^{-1}\kappa_1$ at $(x_0,t_0)$, to apply the tensor maximum principle in Theorem \ref{s2:tensor-mp}, we need to prove the following inequality
\begin{align*}
  Q_1:= &  \alpha F^{-\alpha-2}\ddot{F}^{kl,mn}\nabla_1h_{kl}\nabla_1h_{mn}-\alpha(\alpha+1)F^{-\alpha-3}(\nabla_1F)^2\nonumber\\
&\quad +\alpha(\alpha+1)F^{-\alpha-4}\dot{F}^{kl}\nabla_kF\nabla_lF\kappa_1\\
&\quad+2\alpha F^{-\alpha-1}\sup_{\Lambda}\dot{F}^{kl}\left(2\Lambda_k^p\nabla_lG_{1p}-\Lambda_k^p\Lambda_l^qG_{pq}\right)~\geq~0
\end{align*}
at $(x_0,t_0)$.  Note that  $G_{11}=0$, $\nabla_k G_{11}=0$ at $(x_0,t_0)$, the supremum over $\Lambda$ can be computed exactly as follows:
\begin{align*}
  2\dot{F}^{kl}&\left(2\Lambda_k^p\nabla_lG_{1p}-\Lambda_k^p\Lambda_l^qG_{pq}\right)  \\
  &\quad =2\sum_{k=1}^n\sum_{p=2}^n\dot{f}^k\left(2\Lambda_k^p\nabla_kG_{1p}-(\Lambda_k^p)^2G_{pp}\right)\\
  &\quad =2 \sum_{k=1}^n\sum_{p=2}^n\dot{f}^k\left(\frac{(\nabla_kG_{1p})^2}{G_{pp}}-\left(\Lambda_k^p-\frac{\nabla_kG_{1p}}{G_{pp}}\right)^2G_{pp}\right).
\end{align*}
It follows that the supremum is obtained by choosing $\Lambda_k^p=\frac{\nabla_kG_{1p}}{G_{pp}}$.
Since $\nabla_kG_{1l}=F^{-1}\nabla_kh_{1l}$ for $1\neq l$ at $(x_0,t_0)$, the required inequality for $Q_1$ now becomes:
\begin{align}\label{s6:Q1-1}
  \frac{Q_1}{\alpha F^{-\alpha-2}}=&~ \ddot{F}^{kl,mn}\nabla_1h_{kl}\nabla_1h_{mn}-(\alpha+1)F^{-1}(\nabla_1F)^2\nonumber\\
&\quad +(\alpha+1)F^{-2}\dot{f}^{k}(\nabla_kF)^2\kappa_1 +2 \sum_{k=1}^n\sum_{l=2}^n\frac{\dot{f}^k}{\kappa_l-\kappa_1}(\nabla_kh_{1l})^2~\geq~0.
\end{align}
Using \eqref{s2:F-ddt} to express the second derivatives of $F$, the concavity of $f_*$ and \eqref{s2:inv-conc-1}, the first term of \eqref{s6:Q1-1} can be estimates as follows:
\begin{align*}
\ddot{F}^{kl,pq}\nabla_1h_{kl}\nabla_1h_{pq} =&\ddot{f}^{kl}\nabla_1h_{kk}\nabla_1h_{ll}+2\sum_{k>l}\frac{\dot{f}^k-\dot{f}^l}{\kappa_k-\kappa_l}(\nabla_1h_{kl})^2\\
\geq&~2f^{-1}\dot{f}^k\dot{f}^l\nabla_1h_{kk}\nabla_1h_{ll}-2\frac{\dot{f}^k}{\kappa_k}(\nabla_1h_{kk})^2+2\sum_{k>l}\frac{\dot{f}^k-\dot{f}^l}{\kappa_k-\kappa_l}(\nabla_1h_{kl})^2\\
=&~2F^{-1}(\nabla_1F)^2-2\sum_{k}\frac{\dot{f}^k}{\kappa_k}(\nabla_1h_{kk})^2+2\sum_{k>l}\frac{\dot{f}^k-\dot{f}^l}{\kappa_k-\kappa_l}(\nabla_1h_{kl})^2.
\end{align*}
The fact that $\nabla_kG_{11}=0$ at $(x_0,t_0)$ implies that $\nabla_kF=\kappa_1^{-1}F\nabla_kh_{11}$ at $(x_0,t_0)$. We now have
\begin{align}\label{s6:Q1-2}
 \frac{Q_1}{\alpha F^{-\alpha-2}}\geq & ~(1-\alpha)F^{-1}(\nabla_1F)^2-2\sum_{k}\frac{\dot{f}^k}{\kappa_k}(\nabla_1h_{kk})^2+2\sum_{k>l}\frac{\dot{f}^k-\dot{f}^l}{\kappa_k-\kappa_l}(\nabla_1h_{kl})^2\nonumber\\
&\quad +(\alpha+1)F^{-2}\dot{f}^{k}(\nabla_kF)^2\kappa_1 +2 \sum_{k=1}^n\sum_{l=2}^n\frac{\dot{f}^k}{\kappa_l-\kappa_1}(\nabla_1h_{kl})^2\nonumber\displaybreak[0]\\
=&~(1-\alpha)F\frac{(\nabla_1h_{11})^2}{\kappa_1^2}-2\frac{\dot{f}^1}{\kappa_1}(\nabla_1h_{11})^2-2\sum_{k>1}\frac{\dot{f}^k}{\kappa_k}(\nabla_1h_{kk})^2\nonumber\\
&\quad +2\sum_{k>1}\frac{\dot{f}^k-\dot{f}^1}{\kappa_k-\kappa_1}(\nabla_kh_{11})^2+2\sum_{k>l>1}\frac{\dot{f}^k-\dot{f}^l}{\kappa_k-\kappa_l}(\nabla_1h_{kl})^2\nonumber\\
&\quad +(\alpha+1)\left(\frac{\dot{f}^1}{\kappa_1}(\nabla_1h_{11})^2+\sum_{k>1}\frac{\dot{f}^k}{\kappa_1}(\nabla_kh_{11})^2\right)\nonumber\displaybreak[0]\\ &\quad +2 \left(\sum_{k>1}\frac{\dot{f}^1}{\kappa_k-\kappa_1}(\nabla_kh_{11})^2+\sum_{k>1,l>1}\frac{\dot{f}^k}{\kappa_l-\kappa_1}(\nabla_1h_{kl})^2\right)\nonumber\displaybreak[0]\\
=&~(1-\alpha)\sum_{k>1}\frac{\dot{f}^k\kappa_k}{\kappa_1^2}(\nabla_1h_{11})^2+ \sum_{k>1}\left(2\frac{\dot{f}^k}{\kappa_k-\kappa_1}+(\alpha+1)\frac{\dot{f}^k}{\kappa_1}\right)(\nabla_kh_{11})^2\nonumber\\
&+2\sum_{k>1,l>1}\frac{\dot{f}^k}{\kappa_l-\kappa_1}(\nabla_1h_{kl})^2-2\sum_{k>1}\frac{\dot{f}^k}{\kappa_k}(\nabla_1h_{kk})^2 +2\sum_{k>l>1}\frac{\dot{f}^k-\dot{f}^l}{\kappa_k-\kappa_l}(\nabla_1h_{kl})^2,
\end{align}
where in the second equality we used $F=\sum_{k=1}^n\dot{f}^k\kappa_k$ due to the homogeneity of $F$. Since $f$ is inverse concave, the inequality \eqref{s2:inv-conc} gives that
\begin{equation*}
 2\sum_{k>l>1}\frac{\dot{f}^k-\dot{f}^l}{\kappa_k-\kappa_l}(\nabla_1h_{kl})^2\geq~-2\sum_{k\neq l>1}\frac{\dot{f}^k}{\kappa_l}(\nabla_1h_{kl})^2.
\end{equation*}
There holds
\begin{equation*}
\mathrm{  last ~line ~of ~} \eqref{s6:Q1-2}~\geq~2\sum_{k>1,l>1}\left(\frac{\dot{f}^k}{\kappa_l-\kappa_1}-\frac{\dot{f}^k}{\kappa_l}\right)(\nabla_1h_{kl})^2~\geq~0.
\end{equation*}
Therefore $Q_1\geq 0$ if the power $\alpha\leq 1$.

In summary, if the power $\alpha\in (0,1]$ for $K_N=0,-1$ and $\alpha=1$ for $K_N=1$, we can apply the tensor maximum principle to conclude that the minimum of the smallest eigenvalue of $F^{-1}h_i^j$ over $M_t$ is increasing in time $t$ along the flow \eqref{s1:flow1}. The strictly increasing is due to the strong maximum principle, since otherwise there exists a unit parallel vector field $v$ on $M_{t_0}$ such that $F^{-1}h_{ij}v_iv^i=C$, i.e., $\kappa_1=CF$ on $M_{t_0}$. The constant $C$ must be equal to $1/n$ since if not, $Q_0v^iv^j>0$ on $M_{t_0}$. Then we conclude that $\kappa_n/\kappa_1=1$ everywhere on $M_{t_0}$ and $M_{t_0}$ is a totally geodesic sphere. This completes the proof.
\endproof

\begin{lem}[\cite{Andrews-McCoy-Zheng}]\label{s3:lem2}
If $f_*$ approaches zero on the boundary of $\Gamma_+$, then for any $C>0$ there exists $C'>0$ such that if $\tau\in\Gamma_+$ and $\tau_{\max}\leq Cf_*(\tau)$, then $\tau_{\max}\leq C'\tau_{\min}$.
\end{lem}
By Lemma \ref{s3:lem-pinc}, the principal curvature $\kappa=(\kappa_1,\cdots,\kappa_n)$ of $M_t$ satisfies $ \kappa_{\min}\geq Cf(\kappa)$ along the flow \eqref{s1:flow1}, where $\kappa_{\min}=\min\{\kappa_i\}$ and $C$ is a constant depending only on the initial date $M_0$. Let $\tau_i=1/{\kappa_i}$. Then $\tau=(\tau_1,\cdots,\tau_n)\in \Gamma_+$ and $ \tau_{\max}\leq Cf_*(\tau)$. Applying Lemma \ref{s3:lem2}, we have
 \begin{equation}\label{s3:pinc1}
   \tau_{\max}\leq C\tau_{\min}
 \end{equation}
for all $t\in [0,T)$. Since the principal curvature $\kappa_i=1/{\tau_i}$, we arrive at the following pinching estimate on  the principal curvatures $\kappa=(\kappa_1,\cdots,\kappa_n)$ of $M_t$:
 \begin{equation}\label{s3:pinc2}
   \kappa_{\max}\leq C\kappa_{\min}
 \end{equation}
for all $t\in[0,T)$, where $C$ is a constant depending only on the initial date $M_0$.

\section{Proof of Corollary \ref{s1:thm-conv}}\label{sec:flow-R}
Due to the pinching estimate in the previous section, we can now follow the similar procedure in \cite{Gerhardt2015,Scheuer2015-1} to prove the convergence of  the flow \eqref{s1:flow1}. Since $f$ is inverse concave but may not be concave with respect to its argument, the key $C^{2,\alpha}$ estimate of the flow can not be derived from the usual radial graphical representation. Instead we adopt the Gauss map parametrization of convex hypersurface and write the flow \eqref{s1:flow1} as a parabolic equation of the support function which is concave with respect to the second spatial derivatives due to the inverse concavity of $f$. This idea has been used in \cite{And2000,Andrews-McCoy-Zheng,Urbas91JDG} for curvature flows of convex hypersurfaces in the Euclidean space. In the following, we will describe this procedure in the hyperbolic space and in sphere.

\subsection{The flow in hyperbolic space}\label{sec:4-2}
The convergence of the flow \eqref{s1:flow1} in hyperbolic space follows the proof in \cite{Scheuer2015-1}. The only step that we need to change is to derive the H\"{o}lder estimate on the second derivatives of the solution $M_t$.  As the Euclidean case in \cite{And2000,Andrews-McCoy-Zheng,Urbas91JDG}, we need to derive a parabolic equation of the support function via the Gauss map parametrization which is concave with respect to the second spatial derivatives. Such Gauss map parametrization of curvature flows in hyperbolic space has been formulated recently by the author with Andrews \cite{And-Wei2017} to study the quermassintegral preserving curvature flow in hyperbolic space.

Firstly, we briefly review the Gauss map parametrization of curvature flows in hyperbolic space briefly and refer the readers to \cite{And-Wei2017} for details. Denote by $\mathbb{R}^{1,n+1}$ the Minkowski spacetime, that is the vector space $\mathbb{R}^{n+2}$ endowed with the Minkowski spacetime metric $\langle \cdot,\cdot\rangle$ by
\begin{equation*}
  \langle X,X\rangle~=~-X_0^2+\sum_{i=1}^{n+1}X_i^2
\end{equation*}
for any vector  $X=(X_0,X_1,\cdots,X_{n+1})\in \mathbb{R}^{n+2}$. The hyperbolic space $\mathbb{H}^{n+1}$ is then given by
\begin{equation*}
  \mathbb{H}^{n+1}=~\{X\in \mathbb{R}^{1,n+1},~~\langle X,X\rangle=-1,~X_0>0\}.
\end{equation*}
An embedding $X:M^n\to \mathbb{H}^{n+1}$ of an $n$-dimensional hypersurface induces an embedding $Y:M^n\to B_1(0)\subset \mathbb{R}^{n+1}$ by
\begin{equation}\label{s5:XY}
  X~=~\frac{(1,Y)}{\sqrt{1-|Y|^2}}.
\end{equation}
Let $g_{ij}^X, h_{ij}^X$ and $g_{ij}^Y, h_{ij}^Y$ be the induced metrics and second fundamental forms of $X(M^n)\subset \mathbb{H}^{n+1}$ and $Y(M^n)\subset \mathbb{R}^{n+1}$ respectively, and $N\in \mathbb{R}^{n+1}$ be the unit normal vector of $Y(M^n)$. We have
\begin{align}\label{s5:h-2}
  h_{ij}^X= & \frac{h_{ij}^Y}{\sqrt{(1-|Y|^2)(1-\langle N,Y\rangle^2)}},
\end{align}
\begin{align*}
  g_{ij}^X =& \frac 1 {1-|Y|^2}\left(g_{ij}^Y+\frac{\langle Y,\partial_iY\rangle\langle  Y,\partial_jY\rangle }{(1-|Y|^2)}\right).
\end{align*}

Suppose $X:M^n\times[0,T)\to \mathbb{H}^{n+1}$ is a solution to the flow \eqref{s1:flow1}. Up to a tangential diffeomorphism, the corresponding embedding $Y: M^n\times[0,T)\to \mathbb{R}^{n+1}$ related by \eqref{s5:XY} satisfies the following evolution equation:
\begin{equation}\label{s5:Y-evl}
 \partial_tY=~\sqrt{(1-|Y|^2)(1-\langle N,Y\rangle^2)}F^{-\alpha}(\mathcal{W}^X) N
\end{equation}
where $\mathcal{W}^X$ is the Weingarten matrix of $X(M^n,t)\subset \mathbb{H}^{n+1}$ which has inverse matrix $\mathcal{W}_X^{-1}$
\begin{align}\label{s5:W-inv}
  (\mathcal{W}_X^{-1})_{i}^j =& (h^{-1}_X)^{jk}g_{ki}^X \nonumber\\
   =&(h^{-1}_Y)^{kj}\left(g_{ki}^Y+\frac{\langle  Y,\partial_iY\rangle\langle Y,\partial_kY\rangle }{(1-|Y|^2)}\right)\sqrt{\frac{1-\langle N,Y\rangle^2}{1-|Y|^2}}.
\end{align}

Since each $M_t$ is strictly convex in $\mathbb{H}^{n+1}$, the equation \eqref{s5:h-2} implies that each $Y_t=Y(M^n,t)$ is strictly convex in $\mathbb{R}^{n+1}$ as well.  It's well known that convex hypersurfaces in the Euclidean space can be parametrized via the Gauss map. Given a smooth strictly convex hypersurface $M$ in $\mathbb{R}^{n+1}$, the support function $s: \mathbb{S}^n\to\mathbb{R}$ of $M$ is defined by $s(z)=\sup\{\langle x,z\rangle:x\in\Omega\}$, where $\Omega$ is the convex body enclosed by $M$. Then the hypersurface $M$ is given by the embedding (see \cite{And2000})
\begin{equation*}
  Y(z)=s(z)z+\bar{\nabla}s(z),
\end{equation*}
where $\bar{\nabla}$ is the gradient with respect to the round metric $\bar{g}_{ij}$ on $\mathbb{S}^n$. The principal radii of curvature $\tau=(\tau_1,\cdots,\tau_n)$ are the eigenvalues of
\begin{equation*}
  \tau_{ij}=\bar{\nabla}_i\bar{\nabla}_js+\bar{g}_{ij}s
\end{equation*}
 with respect to $\bar{g}_{ij}$. It can be checked that $\tau_i=1/{\kappa_i}$ for each $i=1,\cdots,n$. Therefore the solution of \eqref{s5:Y-evl} is given up to a tangential diffeomorphism by solving the following scalar parabolic equation on $\mathbb{S}^n$
\begin{equation}\label{s5:Y-evl-Gau}
  \partial_ts=~\sqrt{(1-s^2-|\bar{\nabla}s|^2)(1-s^2)}F_{*}^{\alpha}(\mathcal{W}_X^{-1})
\end{equation}
for the support function $s(z,t)$, where $\mathcal{W}_X^{-1} $ is the matrix \eqref{s5:W-inv} which can be rewritten as
\begin{align*}
  (\mathcal{W}_X^{-1})_{i}^j    =&~\left(\bar{g}^{jq}+\frac{\langle \bar{g}^{ja}\bar{\nabla}_as,\bar{g}^{qb}\bar{\nabla}_bs\rangle}{1-s^2-|\bar{\nabla}s|^2}\right)\tau_{qi}\sqrt{\frac{1-s^2}{1-s^2-|\bar{\nabla}s|^2}}
\end{align*}
in local coordinates on $\mathbb{S}^n$ in terms of $s,\bar{\nabla}s$ and $\tau_{ij}=\bar{\nabla}\bar{\nabla}_js+s\bar{g}_{ij}$.

Since we already have the pinching estimate \eqref{s1:pinc2}, the argument in \cite[\S 3]{Scheuer2015-1} implies that the curvature function $F$ satisfies
\begin{equation*}
  0<c^{-1}\leq F(x,t)\leq c,\qquad \forall~x\in M_t,\quad t\in [0,T)
\end{equation*}
for some constant $c$ depending only on $M_0$ and $\alpha$. This implies the uniform positive two-sides bounds on $\tau_i=1/{\kappa_i}$ and $(\dot{F}_*^{pq})$ for all $t\in [0,T)$. If the maximum existence time $T$ of the flow \eqref{s1:flow1} is finite, comparing with spherical solution implies that the solution $M_t$ stays in a compact subset of $\mathbb{H}^{n+1}$ and $|Y|^2=s^2+|\bar{\nabla}s|^2\leq C<1$ for some constant $C$. Then the equation \eqref{s5:Y-evl-Gau} is uniformly parabolic. Moreover, by the concavity of $F_*$ and $\alpha\leq 1$, the right hand side of the equation \eqref{s5:Y-evl-Gau} is concave with respect to the spatial second derivative $\bar{\nabla}^2s$. Since we have uniform $C^2$ estimate on the support function $s$ in space-time, the H\"{o}lder estimate of Krylov and Evans on second derivatives and the parabolic Schauder theory yield the higher order estimate. A standard continuation argument implies that the flow \eqref{s1:flow1} in hyperbolic space exists for all time $t\in[0,\infty)$ and the evolving hypersurface $M_t$ expands to infinity.

By Lemma \ref{s3:lem-pinc}, the minimum of the smallest eigenvalue of $F^{-1}h_i^j$ over $M_t$ is increasing in time $t$ along the flow \eqref{s1:flow1}. Since the smallest eigenvalue of $F^{-1}h_i^j$ is always smaller than $1/n$ and the flow exists for all time $t\in [0,\infty)$, there exists a constant $\delta\leq 1/n$ such that the minimum of the smallest eigenvalue of $F^{-1}h_i^j$ approaches to $\delta$ as $t\to\infty$. The constant $\delta$ must be equal to $1/n$. If $\delta<1/n$, for any $\delta'\leq \delta$ which is sufficiently close to $\delta$, there exists a sufficiently large time $t_0$ such that minimum of the smallest eigenvalue of $F^{-1}h_i^j$ over $M_{t_0}$ is equal to $\delta'$ and is achieved at $(p,v)\in TM_{t_0}$. The zero order terms \eqref{s3:Q0-1} of the evolution of $G_{ij}$ at $(p,v,t_0)$ satisfy
\begin{align*}
  Q_0v^iv^j\geq &~(\alpha+1)F^{-\alpha-1}\left(F^{-1}\kappa_1\sum_k\dot{f}^k\kappa_k^2-\kappa_1^2\right) \\
  = &~ (\alpha+1)\delta' F^{-\alpha-1}\sum_k\dot{f}^k\kappa_k\left(\kappa_k-\kappa_1\right)\\
  \geq &~ (\alpha+1) F^{-\alpha-1}\dot{f}^n\kappa_n\kappa_1\left(\frac 1{n}-\delta' \right)\\
  \geq &~ (\alpha+1)F^{-\alpha-1}\dot{f}^n\kappa_n\kappa_1\left(\frac 1{n}-\delta\right)~\geq C>0,
\end{align*}
as $\delta<1/n$ and $F,\kappa_i, \dot{f}^i$ are all uniformly bounded. This implies that the minimum of $F^{-1}\kappa_1$ would increase to a constant which is bigger than $1/n$ for a short time, contradicting with the fact $F^{-1}\kappa_1\leq 1/n$. Thus we conclude that $\delta=1/n$ and the pinching ratio $\kappa_n/\kappa_1$ approaches $1$ when $t\to\infty$. The argument in \cite{Scheuer2015-1} (see also \cite{LiWangWei2016}) improves the decay to
\begin{equation*}
  |\kappa_i-1|\leq~c e^{-\frac 2{n^{\alpha}}t},\quad \forall~t>0.
\end{equation*}
The smooth convergence of the flow then follows from the same argument in \cite{Scheuer2015-1}.

\subsection{The flow in sphere}
The pinching estimate \eqref{s1:pinc2} is also the key ingredient to prove the convergence of the flow \eqref{s1:flow1} in sphere. We will employ the dual flow which was introduced by Gerhardt \cite{Gerhardt2015}. Since each $M_t$ is strictly convex in $\mathbb{S}^{n+1}$, considering $M_t$ as a codimension $2$ submanifold in $\mathbb{R}^{n+2}$, the Gauss map $\tilde{X}_t\in T_{X_t}\mathbb{R}^{n+2}$ represents the unit normal vector $\nu\in T_{X_t}\mathbb{S}^{n+1}$ to $M_t$. The mapping $\tilde{X}_t:M^n\to \mathbb{S}^{n+1}$ is also an embedding of a closed strictly convex hypersurface $\tilde{M}_t$. $\tilde{M}_t$ is called the polar set of $M_t$. The Weingarten matrix $\tilde{\mathcal{W}}$ of $\tilde{M}_t$ is the inverse matrix of $\mathcal{W}$ of $M_t$ and the principal curvatures $\tilde{\kappa}=(\tilde{\kappa}_1,\cdots,\tilde{\kappa}_n)$ of $\tilde{M}_t$ satisfy $\tilde{\kappa}_i=1/{\kappa_i}$, where $\kappa_i$ are principal curvatures of $M_t$.

In \cite{Gerhardt2015}, Gerhardt proved that the polar set $\tilde{M}_t$ of $M_t$ satisfies the following contracting curvature flow
\begin{equation}\label{s5:dual flow}
 \frac{\partial}{\partial t}\tilde{X}=~-F_*(\mathcal{\tilde{W}})\tilde{\nu}
 \end{equation}
in $\mathbb{S}^{n+1}$. The flow hypersurfaces of \eqref{s1:flow1} and \eqref{s5:dual flow} are polar sets of each other. Since the principal curvatures $\tilde{\kappa}=(\tilde{\kappa}_1,\cdots,\tilde{\kappa}_n)$ of $\tilde{M}_t$ satisfy $\tilde{\kappa}_i=1/{\kappa_i}$, the pinching estimate on $\kappa_i$ of $M_t$ yields the pinching estimate on $\tilde{\kappa}_i$ of $\tilde{M}_t$. Moreover, since $F(\mathcal{W})$ is inverse concave, the operator $F_*(\tilde{\mathcal{W}})$ is concave with respect to the components of $\tilde{\mathcal{W}}$. Hence the H\"{o}lder estimate of Krylov and Evans on the second derivatives can be applied. Then the same argument as in \cite{Gerhardt2015} implies that the flow \eqref{s5:dual flow} contracts to a point in finite time and properly rescaled hypersurfaces have uniformly $C^k$ estimate for all $k\geq 2$. Moreover, if $F_*$ is strictly concave or $F_*$ equals to the mean curvature $H$, the properly rescaled flow hypersurfaces converge to a geodesic sphere smoothly.  Back to our expanding flow \eqref{s1:flow1} in sphere, we have that the flow expands to the equator smoothly in finite time, and in the case that $F_*$ is strictly concave or $F_*$ equals to the mean curvature $H$, the properly rescaled hypersurfaces converge to a geodesic sphere continuously.

To show the smooth convergence of the flow \eqref{s1:flow1}, we still need to derive higher order estimates on the solution $M_t$. This can be done similarly as in \S \ref{sec:4-2} via a Gauss map parametrization. The embedding  $X:M^n\to \mathbb{S}^{n+1}\subset \mathbb{R}^{n+2}$ induces an embedding $Y:M^n\to \mathbb{R}^{n+1}$ by
\begin{equation}\label{s6:XY}
  X~=~\frac{(1,Y)}{\sqrt{1+|Y|^2}}.
\end{equation}
Suppose that $X(M^n)$ is strictly convex, then $X(M^n)$ lies in a hemisphere and the corresponding $Y(M^n)$ is a hypersurface in $\mathbb{R}^{n+1}$ which is also strictly convex since
 \begin{align}\label{s6:h-2}
  h_{ij}^X= & \frac{h_{ij}^Y}{\sqrt{(1+|Y|^2)(1+\langle N,Y\rangle^2)}},
\end{align}
where $h_{ij}^X$ and $h_{ij}^Y$ are the second fundamental forms of $X(M^n)\subset \mathbb{S}^{n+1}$ and $Y(M^n)\subset \mathbb{R}^{n+1}$ respectively, and $N\in \mathbb{R}^{n+1}$ is the unit normal of $Y(M^n)$. Suppose $X:M^n\times[0,T)\to \mathbb{S}^{n+1}$ is a solution to the flow \eqref{s1:flow1}. Up to a tangential diffeomorphism, the corresponding embedding $Y: M^n\times[0,T)\to \mathbb{R}^{n+1}$ related by \eqref{s6:XY} satisfies the following evolution equation:
\begin{equation}\label{s6:Y-evl}
 \partial_tY=~\sqrt{(1+|Y|^2)(1+\langle N,Y\rangle^2)}F^{-1}(\mathcal{W}^X) N,
\end{equation}
where $\mathcal{W}^X$ is the Weingarten matrix of $X(M^n,t)\subset \mathbb{S}^{n+1}$. Since each $M_t$ is strictly convex in $\mathbb{S}^{n+1}$, the equation \eqref{s6:h-2} implies that each $Y_t=Y(M^n,t)$ is strictly convex in $\mathbb{R}^{n+1}$ as well.  As in \S \ref{sec:4-2}, the solution of \eqref{s6:Y-evl} is then given up to a tangential diffeomorphism by solving the following scalar parabolic equation on $\mathbb{S}^n$
\begin{equation}\label{s6:Y-evl-Gau}
  \partial_ts=~\sqrt{(1+s^2+|\bar{\nabla}s|^2)(1+s^2)}F_{*}(\mathcal{W}_X^{-1})
\end{equation}
for the support function $s(z,t)$, where $\mathcal{W}_X^{-1}$ is the inverse matrix of $\mathcal{W}^X$ and is given by
\begin{align*}
  (\mathcal{W}_X^{-1})_{ij}    =&~\left(\bar{g}^{jq}-\frac{\langle \bar{g}^{ja}\bar{\nabla}_as,\bar{g}^{qb}\bar{\nabla}_bs\rangle}{1+s^2+|\bar{\nabla}s|^2}\right)\tau_{qi}\sqrt{\frac{1+s^2}{1+s^2+|\bar{\nabla}s|^2}}
\end{align*}
in local coordinates on $\mathbb{S}^n$ in terms of $s,\bar{\nabla}s$ and $\tau_{ij}=\bar{\nabla}_i\bar{\nabla}_js+s\bar{g}_{ij}$. The equation \eqref{s6:Y-evl-Gau} is parabolic and is concave with respect to the second spatial derivatives. The argument in \cite{Gerhardt2015} can be used to deduce the required higher regularity.



\begin{thebibliography}{99}
\bibitem{And1994-2}
Ben Andrews, \emph{Contraction of convex hypersurfaces in {E}uclidean space},
  Calc. Var. Partial Differ. Equ. \textbf{2} (1994), no.~2, 151--171.

\bibitem{And1994-3}
Ben Andrews, \emph{Contraction of convex hypersurfaces in Riemannian spaces}, J.
  Differential Geom. \textbf{39} (1994), no.~2, 407--431.

\bibitem{And2000}
Ben Andrews, \emph{Motion of hypersurfaces by Gauss curvature,} Pacific. J. Math. \textbf{195}(1), 1--34 (2000)

\bibitem{And04}
Ben Andrews,  \emph{Fully nonlinear parabolic equations in two space variables},
  available at arXiv: math.DG/0402235 (2004).

\bibitem{Andrews2007}
Ben Andrews,  \emph{Pinching estimates and motion of hypersurfaces by curvature
  functions}, J. Reine Angew. Math. \textbf{608} (2007), 17--33.

\bibitem{ALM14}Ben Andrews, Mat Langford and James McCoy, \emph{Convexity estimates for hypersurfaces
moving by convex curvature functions}, Analysis and PDE, \textbf{7 }(2014), No. 2, 407--433.

\bibitem{Andrews-McCoy-Zheng}
Ben Andrews, James McCoy, and Yu~Zheng, \emph{Contracting convex hypersurfaces
  by curvature}, Calc. Var. Partial Differ. Equ. \textbf{47} (2013), no.~3-4, 611--665.

\bibitem{And-Wei2017}
Ben Andrews and Yong Wei,  \emph{Quermassintegral preserving curvature flow in Hyperbolic space}, arXiv:1708.09583.

\bibitem{Chow-Gu96}B. Chow and R. Gulliver,  \emph{Aleksandrov reflection and nonlinear evolution equations. I. The $n$-sphere and $n$-ball}, Calc. Var. Partial Differ. Equ. \textbf{4}(3), 249--264 (1996)

\bibitem{Ding-11}Qi Ding, \emph{The inverse mean curvature flow in rotationally symmetric spaces}. Chin. Ann. Math.
Ser. B \textbf{32} (2011), no. 1, 27--44.

\bibitem{Gerhardt1990}
Claus Gerhardt, \emph{Flow of nonconvex hypersurfaces into spheres}, J.
  Differential Geom. \textbf{32} (1990), no.~1, 299--314.

\bibitem{Gerhardt2006}
Claus Gerhardt, \emph{Curvature problems}, Series in Geometry and Topology, vol.~39,
  International Press, Somerville, MA, 2006.


\bibitem{gerhardt2011-H^n}
Claus Gerhardt, \emph{Inverse curvature flows in hyperbolic space}, J. Differential
  Geom. \textbf{89} (2011), no.~3, 487--527.

\bibitem{Gerhardt2014}
Claus Gerhardt, \emph{Non-scale-invariant inverse curvature flows in {E}uclidean
  space}, Calc. Var. Partial Differ. Equ. \textbf{49} (2014),
  no.~1-2, 471--489.

\bibitem{Gerhardt2015}
Claus Gerhardt, \emph{Curvature flows in the sphere}, J. Differential Geom.
  \textbf{100} (2015), no.~2, 301--347.

\bibitem{Ha1982}
R. S. Hamilton,\emph{ Three-manifolds with positive Ricci curvature}, J. Differential. Geom. \textbf{17}(2), 255--306 (1982)

\bibitem{Krylov} N.V. Krylov, \emph{Boundedly inhomogeneous elliptic and parabolic equations}. Izv. Akad. Nauk SSSR Ser.
Mat. 46(3), 487--523, 670 (1982) (Russian)


\bibitem{Kron-Sche2017}Heiko Kr\"{o}ner and Julian Scheuer, \emph{Expansion of pinched hypersurfaces of the Euclidean and hyperbolic space by high powers of curvature}, arXiv:1703.07087.

\bibitem{LiWangWei2016}
Haizhong Li, Xianfeng Wang and Yong Wei, \emph{Surfaces expanding by non-concave curvature functions}, arXiv:1609.00570.
%

\bibitem{Q-Li2010}
Qi-Rui Li, \emph{Surfaces expanding by the power of the {G}auss curvature
  flow}, Proc. Amer. Math. Soc. \textbf{138} (2010), no.~11, 4089--4102.

\bibitem{lieberman1996second}
Gary~M Lieberman, \emph{Second order parabolic differential equations}, World
  scientific, 1996.

\bibitem{Mcc2017}James A. McCoy, \emph{More mixed volume preserving curvature flows}, J. Geom. Anal. (online first), 2017.

\bibitem{makowski2013rigidity}
Matthias Makowski and Julian Scheuer, \emph{Rigidity results, inverse curvature
  flows and alexandrov-fenchel type inequalities in the sphere},
  Asian J. Math., \textbf{20}, no. 5, p. 869--892, (2016).


\bibitem{scheuer2015-gradient}
Julian Scheuer, \emph{Gradient estimates for inverse curvature flows in
  hyperbolic space}, Geometric Flows \textbf{1} (2015), no.~1, 11--16.

\bibitem{Scheuer2015-1}
Julian Scheuer, \emph{Non-scale-invariant inverse curvature flows in hyperbolic
  space}, Calc. Var. Partial Differ. Equ. \textbf{53} (2015),
  no.~1-2, 91--123.

\bibitem{Sch2016} Julian Scheuer, \emph{Pinching and asymptotical roundness for inverse curvature flows in Euclidean space,} J. Geom. Anal. \textbf{26} (2016), no. 3, 2265--2281.

\bibitem{Schnurer2006}
Oliver~C. Schn{\"u}rer, \emph{Surfaces expanding by the inverse {G}au\ss\
  curvature flow}, J. Reine Angew. Math. \textbf{600} (2006), 117--134.

\bibitem{urbas1990expansion}
John~IE Urbas, \emph{On the expansion of starshaped hypersurfaces by symmetric
  functions of their principal curvatures}, Math. Z.  \textbf{205} (1990), no.~1, 355--372.

\bibitem{Urbas91JDG}
John~IE Urbas, \emph{An expansion of convex hypersurfaces}, J. Differential Geom., 1991, 33(1): 91--125.
\end{thebibliography}
\end{document}